\documentclass[12pt]{article}
\usepackage{amsmath, amsfonts}

\usepackage{color}
\numberwithin{equation}{section}
\def \no { \noindent}

\def\ds{\displaystyle}

\def\B1{B_{1/2}}

\def\Bl1{{\lambda_1}}
\def\Box{\hfill\rule{2.5mm}{2.5mm}}
\def\C{{\cal {C}}}

\def\H{{\cal H}}

\def\L{{\cal L}}

\def\N{{\mathbb {N}}}

\def\R{{\mathbb {R}}}

\def\RR{{\cal R}}
\def\SS{{\cal S}}

\def\argth{\mathop {\arg \tanh}}

\def\bl1{{\bar \lambda_1}}

\def\build#1_#2^#3{\mathrel{
\mathop{\kern 0pt#1}\limits_{#2}^{#3}}}

\def\d{\displaystyle}

\def\E{{\cal E}}

\def\h1{\mathop{\rm H^1_{\rm loc,\rm u}}}

\def\iint{\displaystyle\int_{-1}^1}

\def\l2{\mathop{\rm L^2_{\rm loc,\rm u}}}

\def\pr{{\partial_r}}
\def\ps{\partial_s}

\def\py{\partial_y}

\newcommand{\aref}[1]{(\ref{#1})}

\newcommand{\vc}[2]{
\left(
\begin{array}{l}
#1\\
#2
\end{array}
\right)
}

\newtheorem{cor}{Corollary}[section]

\newtheorem{lem}[cor]{Lemma}
\newtheorem{prop}[cor]{Proposition}

\newtheorem{propo}{Proposition}
\newtheorem{theor}[propo]{Theorem}
\newtheorem{coro}[propo]{Corollary}

\begin{document}

\title{\bf Blow-up behavior  for the Klein-Gordon and other perturbed semilinear wave equations}

\author{M.A. Hamza\\
{\it \small  {Facult\'e des Sciences de Tunis}}\\
H. Zaag \\
 {\it \small   {CNRS UMR 7539 LAGA  Universit\'e Paris 13}} }

\maketitle

\begin{abstract}
We give blow-up results for the Klein-Gordon equation and other perturbations of the semilinear wave equations with superlinear power nonlinearity, in one space dimension or in higher dimension under radial symmetry outside the origin.
\end{abstract}

\medskip


\medskip


\section{Introduction}
We consider one dimensional solutions and higher dimensional  radial
solutions of the following  Klein-Gordon equation:
\begin{equation}\label{equrn}
\left\{
\begin{array}{l}
\partial_t^2 U =\Delta U+|U|^{p-1}U-U,\\
U(0)=U_0\mbox{ and }U_t(0)=U_1,
\end{array}
\right.
\end{equation}
where $U(t):x\in\R^N \rightarrow U(x,t)\in\R$, $U_0\in \rm H^1_{\rm loc,u}$
and $U_1\in \rm L^2_{\rm loc,u}$.\\
 The space $\rm L^2_{\rm loc,u}$ is the set of all $v$ in
$\rm L^2_{\rm loc}$ such that
\[
\|v\|_{\rm L^2_{\rm loc,u}}\equiv\d\sup_{a\in
\R^N}\left(\int_{|x-a|<1}|v(x)|^2{\mathrm{d}}x\right)^{1/2}<+\infty,
\]
 and the space ${\rm H}^1_{\rm loc,u}= \{ v\;|\;v, \nabla v \in {\rm L}^2_{\rm loc,u}\}$.

\bigskip

The  nonlinear Klein Gordon equation appears as a model of self-focusing waves in nonlinear optics (see Bizo\'n, Chamj and Szpak \cite{BCS11}).

\bigskip

More generally, we consider  the following semilinear  wave
equation:
\begin{equation}\label{gen}
\left\{
\begin{array}{l}
\partial_t^2 U =\Delta U+|U|^{p-1}U+f(U)+g(|x|,t,\nabla U.\frac{x}{|x|},\partial_t U ),\\
U(0)=U_0\mbox{ and }U_t(0)=U_1.
\end{array}
\right.
\end{equation}
We assume that the functions $f$  and $g$  are ${\cal {C}}^1$
functions, where $f:\R\rightarrow \R $  and $g:\R^{4}\rightarrow \R
$ satisfy the following conditions
\begin{eqnarray*}
(H_f)& |{f(u)}|\le M(1+|u|^q), \  &{\textrm {for all }}\ y\in \R \  \qquad{{\textrm {with}}}\ \ (q<p,\ \ M>0),\\
(H_g)& |{g(X,t,v,z)}|\le M(1+|v|+|z|), & {\textrm {for all
}}\ X, t,v,z\in \R \  \qquad{{\textrm {with}}}\  (M>0).
\end{eqnarray*}
We assume also  that the function $g$ is globally Lipschitz. Finally,  we assume  that
\begin{equation}\label{condp}
1<p\mbox{ and }p\le 1+\frac 4{N-1}\mbox{ if } N\ge 2.
\end{equation}
Since $U$ is radial if $N\ge2$, we introduce
\begin{equation}\label{defu}
u(r,t) = U(r,t)\mbox { for }r\in \R, \mbox { if }  N=1
\end{equation}
and
\begin{equation}\label{defu1}
u(r,t) = U(x,t)\mbox { if }r=|x| \mbox { and  } N\ge2
\end{equation}
and rewrite \eqref{gen} as
\begin{equation}\label{equ}
\left\{
\begin{array}{l}
\partial^2_{t} u =\partial^2_r u+\frac{(N-1)}r \pr u+|u|^{p-1}u+f(u)+g(r,t,\partial_r u, \partial_{t}u),\\
\pr u(0,t)=0 \mbox { if } N\ge2,\\
u(r,0)=u_0(r)\mbox{ and }u_t(r,0)=u_1(r),
\end{array}
\right.
\end{equation}
where $u(t):r\in I \rightarrow u(r,t)\in\R$, with $I=\R^+$ if $N\ge2$ and $I=\R$ if $N=1$.

\bigskip

The  Cauchy problem  of  equation \eqref{gen} is  solved in
$H^{1}_{loc,u}\times L^{2}_{loc,u}$. This follows from the finite
speed of propagation and the  wellposedness  in $H^{1} \times
L^{2}$ (see  for example Georgiev  and Todorova
 \cite{GTjde94}) , valid whenever $ 1< p<1+\frac{4}{N-2}$. The existence of
blow-up   solutions  $u(t)$ of  \eqref{gen}  follows from energy  techniques (see  for example
Levine  and  Todorova
\cite{LT})
and  Todorova
\cite{Tn00}).

\bigskip

If $u$ is a blow-up solution of \eqref{equ}, we define (see for
example Alinhac \cite{Apndeta95}) a 1-Lipschitz curve
$\Gamma=\{(r,T(r))\}$ where $r\in I $
such that the maximal influence domain $D$ of $u$ (or the domain of definition of $u$) is written as
\begin{equation}\label{defdu}
D=\{(r,t)\;|\; r\in I,\ t< T(r)\}.
\end{equation}
$\Gamma$ is called the blow-up graph of $u$.
A point $r_0\ge 0$ is a non-characteristic point if there are
\begin{equation}\label{nonchar}
\delta_0\in(0,1)\mbox{ and }t_0<T(r_0)\mbox{ such that }
u\;\;\mbox{is defined on }{\cal C}_{r_0, T(r_0), \delta_0}\cap \{t\ge t_0\}\cap\{r\in I\}
\end{equation}
where ${\cal C}_{\bar r, \bar t, \bar \delta}=\{(r,t)\;|\; t< \bar t-\bar \delta|r-\bar r|\}$.
We denote by $\RR\subset I$ (resp. $\SS\subset I$) the set
of non-characteristic (resp. characteristic) points.

\bigskip

In the case $(f,g)\equiv (0,0)$, equation (\ref{gen}) reduces to the semilinear wave equation:
\begin{equation}\label{par}
\left\{
\begin{array}{l}
\partial_t^2 U =\Delta U+|U|^{p-1}U,\\
U(0)=U_0\mbox{ and }U_t(0)=U_1.
\end{array}
\right.
\end{equation}
  In a series of papers \cite{MZjfa07}, \cite{MZcmp08}, \cite{MZajm10} and \cite{MZisol10}
  (see also the note \cite{MZxedp10}), Merle and Zaag  give a full picture of the blow-up
  for solutions of \eqref{par} in one space dimension. Recently, in
\cite{CZarxiv}, C\^ote and Zaag refine some of those results and construct a blow-up solution with a
characteristic point $a$, such that the asymptotic behavior of the
solution near $(a, T (a))$ shows a decoupled sum of $k$ solitons
with alternate signs.
  Moreover, in  \cite{MZbsm11},
     Merle and Zaag extend  all   their results to higher dimensions in the radial case, outside the origin.
Our aim in this work is to generalize the result obtained for
equation   (\ref{par}) in  \cite{MZbsm11} to equation (\ref{gen}).
Let us note  that all our results and proofs hold    in both  cases of (\ref{equ}) ($N=1$ and $N\ge 2$).
 However,  the situation is a bit more delicate when  $N\ge 2$,
 since we have to
 avoid the origin which brings a singular  term
 $\frac{N-1}r\partial_ru$ to (\ref{equ}). Thus, for completeness, we focus  on the case $N\ge 2$ avoiding  $r=0$, and stress the fact that all our results hold in the case    $N=1$, even when $r=0$, and with no symmetry assumptions.

\bigskip

Throughout this paper, we consider $U(x,t)$ a radial blow-up
solution of equation \eqref{gen}, and use the notation $u(r,t)$
introduced in \eqref{defu}.
We proceed in 3 sections:\\
- in Section \ref{seclyap}, we give a new Lyapunov functional for equation \eqref{equ} and bound the solution in the energy space.\\
- in Section \ref{secnonchar}, we study $\RR$, in particular the blow-up behavior of the solution and the regularity of the blow-up set there.\\
- in Section \ref{secchar}, we focus on $\SS$,  from the point
of view of the blow-up behavior, the regularity of the blow-up set
and the construction of  a multi-soliton solution.

\bigskip

 We are aware that our analysis is a
generalization of the radial case of equation \eqref{par} treated by
Merle and Zaag in \cite{MZbsm11}. For that reason, we will give the
statements of the results for equation \eqref{gen} and focus only on
 how to deal with the new perturbation terms appearing in
\eqref{gen}. Let  us add that, we believe that our contribution is
non trivial and introduces a new approach for perturbed problems.
Moreover, it proves a bunch of results, especially for the
Klein-Gordon equation (\ref{equrn}).

\bigskip


\section{A  Lyapunov functional  for equation \eqref{eqw}  and a blow-up criterion in the radial case}\label{seclyap}

\bigskip

We showed in \cite{HZlyap10} and \cite{HZlyapc10} that the argument
of Antonini, Merle and Zaag in \cite{AMimrn01}, \cite{MZajm03},
\cite{MZma05} and \cite{MZimrn05} extends through a perturbation
method to equation \eqref {gen}
 with no gradient terms   even for non radial solutions. The key idea is to modify the Lyapunov functional of
\cite{AMimrn01} with exponentially small terms and define a new
functional which is  decreasing in time and gives a
blow-up criterion. In  \cite{MZbsm11},  Merle and Zaag successfully
used our  ideas to derive  a Lyapunov  functional for the radial
case with no perturbations (i.e. for equation (\ref{equ}) with
$(f,g)\equiv(0,0)$). Here, we further   refine our argument in
\cite{HZlyap10} and \cite{HZlyapc10} to derive a
 Lyapunov  functional for equations bearing the two features: the
 presence of perturbation terms and the radial symmetry.
For the reader's convenience,
 we briefly recall the argument  in the following.

\bigskip

Given $r_0>0$, we recall the following similarity variables' transformation
\begin{equation}\label{defw}
w_{r_0}(y,s)=(T(r_0)-t)^{\frac 2{p-1}}u(r,t),\;\;y=\frac{r-r_0}{T(r_0)-t},\;\;
s=-\log(T(r_0)-t).
\end{equation}
The function $w=w_{r_0}$ satisfies the following equation for all
$y\in (-1,1)$

and $s\ge \max\left(-\log T(r_0), -\log r_0\right)$:
\begin{eqnarray}
\label{eqw}
\partial^2_{s}w&=& \L w-\frac{2(p+1)}{(p-1)^2}w+|w|^{p-1}w
-\frac{p+3}{p-1}\partial_sw-2y\partial^2_{y,s}
w\nonumber\\
&&+e^{-s}\frac{(N-1)}{r_0+ye^{-s}} \py w
+e^{-\frac{2ps}{p-1}}f\Big(e^{\frac{2s}{p-1}}w\Big)\\ &&+
e^{-\frac{2ps}{p-1}}g\Big(r_0+ye^{-s},T_0-e^{-s},e^{\frac{(p+1)s}{p-1}}\partial_yw,e^{\frac{(p+1)s}{p-1}}(\partial_sw+y\partial_y
w+\frac{2}{p-1}w)\Big) ,\nonumber
\end{eqnarray}
\begin{equation}\label{defro}
\mbox{where }\L w = \frac 1\rho \py \left(\rho(1-y^2) \py w\right)\mbox{ and }
\rho(y)=(1-y^2)^{\frac 2{p-1}}.
\end{equation}
In the whole paper, we denote
\begin{equation}\label{F}
F(u)=\int_0^uf(v){\mathrm{d}}v.
\end{equation}

\bigskip

Let us recall  that  for the case $(f,g)\equiv(0,0)$,  the
Lyapunov functional in one space dimension is
\begin{equation}\label{defE}
E_0(w(s))=\iint \left(\frac 12 (\ps w)^2 + \frac 12
\left(\partial_y w\right)^2 (1-y^2)+\frac{(p+1)}{(p-1)^2}w^2 - \frac
1{p+1} |w|^{p+1}\right)\!\rho {\mathrm{d}}y,
\end{equation}
which is defined in the Hilbert space
\begin{equation}\label{defnh0}
\H = \left\{q \in {\rm H^1_{loc}}\times {\rm L^2_{loc}(-1,1)}
\;\;|\;\;\|q\|_{\H}^2\equiv \int_{-1}^1 \left(q_1^2+\left(q_1'\right)^2
 (1-y^2)+q_2^2\right)\rho {\mathrm{d}}y<+\infty\right\}.
\end{equation}
Introducing
\begin{eqnarray}
\label{defF} E(w(s),s)&=&E_0(w(s))+I(w(s),s)+J(w(s),s)
\end{eqnarray}
where,
\begin{eqnarray}
\label{f2} \quad
  I(w(s),s)&=&- e^{-\frac{2(p+1)s}{p-1}}\displaystyle\iint F(e^{\frac{2s}{p-1}}w)\rho
  {\mathrm{d}}y,
\end{eqnarray}
\begin{eqnarray}
\label{f3} \quad
 J(w(s),s)&=& -e^{-\gamma s}\displaystyle\iint w\ps w\rho
 {\mathrm{d}}y,
\end{eqnarray}
 with \begin{eqnarray}\label{gamma} \gamma=\min(\frac12,\frac{p-q}{p-1}
 )>0,
 \end{eqnarray}
we claim the following:
\begin{prop}[A new functional for equation \eqref{eqw}]\label{prophamza}
$ $\\
(i) There exists $C=C(p,N,M,q)>0$ and $S_0(p,N,M,q)\in\R$ such that
for all $r_0>0$ and for all $s\ge \max\left(-\log T(r_0), S_0,
-4\log r_0,-\log \frac{r_0}2\right)$,
\begin{equation}\label{edoF}
\frac d{ds}E(w_{r_0}(s),s) \le \frac{p+3}2 e^{-\gamma
s}E(w_{r_0}(s),s) -\frac 3{p-1}\iint (\ps w_{r_0})^2
\frac{\rho}{1-y^2} {\mathrm{d}}y+Ce^{-2\gamma s}.
\end{equation}
(ii) {\bf (A blow-up criterion)}
 There exists $S_1(p,N,M,q)\in \R$ such that, for all $s\ge
\max\left(s_0, S_1,\right)$, we have $H(w(s),s)\ge 0$.

\end{prop}
{\bf Remark}: From (i), we see that the Lyapunov functional
for equation \eqref{eqw} is in fact $H(w_{r_0}(s),s)$ where
\begin{equation}\label{defH}
 H(w(s),s)=E(w(s),s)e^{\frac{ p+3}{2\gamma} e^{-\gamma s}}+\mu
e^{-2\gamma s}
\end{equation}
 not $E(w_{r_0}(s),s)$ nor $E_0(w_{r_0}(s))$, for some large constant $\mu$. \\
{\bf Remark}:
We already know from \cite{HZlyap10} and \cite{HZlyapc10} that even
in the non-radial setting, equation \eqref{equ} has a Lyapunov
functional given by a perturbed  form of a natural extension to
higher dimensions of $E(w_{r_0}(s),s)$ \eqref{defE}. Unfortunately,
already for the non-perturbed case of (\ref{gen}) with $(f,g)\equiv
(0,0)$, due to the lack of information
 on stationary solutions in similarity variables in dimensions $N\ge 2$,
it wasn't possible to go further in the analysis, and the
investigation
  had to stop at the step of bounding the solution
  in similarity variables. On the contrary, when $N=1$,
  Merle and Zaag  could obtain a very precise characterization of blow-up in \cite{MZjfa07}, \cite{MZcmp08},
  \cite{MZajm10}, \cite{MZisol10} (with some refinements by   C\^ote and  Zaag in \cite{CZarxiv}).\\
 Here, considering perturbations as stated in
(\ref{gen}) and
 restricting ourselves
to one dimensional solutions or higher dimensional
 radial solutions,
   we find a different Lyapunov
 functional.
Considering  arbitrary blow-up points in one-space
dimension(including the origin) and any non-zero blow-up
 in higher dimensions,
  the characterization of stationary solutions in one space dimension is enough, and
 we are able to go in our analysis as far as in the one-dimensional case.

\bigskip

Following  \cite{MZajm03} and \cite{MZjfa07}, together with our
techniques to handle perturbations in  \cite{HZlyap10} and
\cite{HZlyapc10}, we derive with no difficulty the following:
\begin{prop}\label{boundedness}{\bf (Boundedness of the solutions of equation
\eqref{eqw} in the energy space)} For all $r_0>0$,
there is a $C_2(r_0)>0$ and $S_2(r_0)\in\R$ such that for all $r\in[\frac{r_0}2, \frac{3r_0}2]$ and $s\ge S_2(r_0)$,
\[\iint\left((\py w_r)^2 (1-y^2) +  (w_r)^2+ (\partial_s w_r)^2+|w_r|^{p+1}\right) \rho {\mathrm{d}}y\le C_2(r_0).\]
 \end{prop}
{\it Proof}: The adaptation is straightforward from  \cite{MZajm03}
and Proposition 3.5 page 66 in \cite{MZjfa07}. The only difference
is in
 the justification of the limit at infinity of $E_0(w_{r_0}(s))$,
 which follows from the limit of $H(w_{r_0}(s),s)$ defined in \eqref{defH}.
 In fact, we know from Proposition \ref{prophamza} that $H(w_{r_0}(s),s)$ is
 decreasing and bounded from below, and such an information is unavailable for $E_0(w_{r_0}(s))$.
\Box

\bigskip

{\it Proof of Proposition \ref{prophamza}}:\\
(i) Consider $r_0>0$, $s\ge \max(-\log T(r_0),0, -\log \frac{r_0}2,
 -4\log {r_0}
)$ and write $w=w_{r_0}$ for simplicity. From the similarity
variables' transformation \eqref{defw}, we see that
\begin{equation}\label{range}
r=r_0+ye^{-s}\in\left[\frac{r_0}2, \frac{3r_0}2\right].
\end{equation}
Multiplying equation \eqref{eqw} by $\ps w\rho$ and integrating
for $y\in(-1,1)$, we see by  \eqref{defE}  and \eqref{f2}  that
\begin{multline}\label{Eprime}
\frac d{ds}\big(E_0(w(s))+I(w(s),s)\big) =\frac {-4}{p-1} \iint
 \frac{(\ps w)^2} {1-y^2}\rho {\mathrm{d}}y +\underbrace{(N-1)
e^{-s} \iint \ps w\py w
\frac{\rho}r {\mathrm{d}}y}_{I_1(s)}\nonumber\\
+\underbrace{\frac{2(p+1)}{p-1}e^{-\frac{2(p+1)s}{p-1}}\iint
F\Big(e^{\frac{2s}{p-1}}w\Big)\rho {\mathrm{d}}y}_{I_2(s)}
+\underbrace{\frac{2}{p-1}e^{-\frac{2ps}{p-1}}\iint f\Big(e^{\frac{2s}{p-1}}w\Big)w\rho {\mathrm{d}}y }_{I_3(s)}\\
+\underbrace{ e^{-\frac{2ps}{p-1}}\! \iint\!\!
g\Big(r_0+ye^{-s},T_0-e^{-s},e^{\frac{(p+1)s}{p-1}}\partial_yw,e^{\frac{(p+1)s}{p-1}}(\partial_sw+y\partial_y
w+\frac{2w}{p-1})\Big) \ps w\rho
{\mathrm{d}}y}_{I_4(s)}.\nonumber
\end{multline}
where $r$ is defined in \eqref{range}.
Using \eqref{range}, we write
\begin{eqnarray}
\label{CS1} |I_1(s)|\le Ce^{-s}\iint (\py w)^2\rho(1-y^2)
{\mathrm{d}}y + \frac{Ce^{-s}}{r_0^2}\iint (\ps w)^2
\frac{\rho}{1-y^2} {\mathrm{d}}y.
\end{eqnarray}
Using the fact that
\begin{equation}\label{born} |{F(x)}|+|x{f(x)}|\le C(
1+|x|^{q+1})\le C( 1+|x|^{p+1}),
\end{equation}
 where $F$ and $f$ are defined in
\eqref{F} and \eqref{gen}, we obtain that
\begin{eqnarray}\label{I10}
|I_2(s)|+|I_3(s)| &\le& Ce^{-\frac{2(p-q)s}{p-1}}+
Ce^{-\frac{2(p-q)s}{p-1}}\iint |w|^{p+1}\rho {\mathrm{d}}y.
\end{eqnarray}
Using the inequality $ab\le \frac{a^2}{2}+\frac{b^2}{2}$ and the
hypothesis  $ (H_g)$, we write  that
\begin{eqnarray}\label{I122}
|I_4(s)|\le Ce^{-s} \iint \Big((\ps w)^2+(w^2)\Big)\rho
{\mathrm{d}}y+ Ce^{-s}\iint|\py w ||\ps w|\rho {\mathrm{d}}y+
C e^{-s}.\qquad
\end{eqnarray}
Similarly, we prove that
 \begin{eqnarray}\label{I13}
\iint|\py w ||\ps w|\rho {\mathrm{d}}y\le \iint (\ps
w)^2\frac{\rho}{1-y^2}{\mathrm{d}}y +\iint (\py
w)^2(1-y^2)\rho{\mathrm{d}}y.
\end{eqnarray}
Combining  (\ref{I122}) and (\ref{I13}), we conclude that
\begin{eqnarray}\label{I16}
|I_4(s)|\le Ce^{-s} \iint\!\! \Big( (\py w)^2(1-|y|^2) +
\frac{(\ps w)^2}{1-y^2}+w^2\Big)\rho {\mathrm{d}}y+
Ce^{-s}.
\end{eqnarray}
Then, by using  (\ref{Eprime}), (\ref{CS1}), (\ref{I10}) and
(\ref{I16}), we deduce that
\begin{eqnarray}\label{E}
\frac d{ds}\big(E_0(w(s))+I(w(s),s)\big) &\le&(-\frac
4{p-1}+Ce^{-\frac{s}2})
\iint (\ps w)^2 \frac{\rho}{1-y^2} {\mathrm{d}}y\nonumber\\
&& + Ce^{-s}\iint
\Big( (\py w)^2(1-|y|^2) +w^2\Big)\rho {\mathrm{d}}y \nonumber\\
&&+ Ce^{-2 \gamma s}\iint |w|^{p+1}\rho
{\mathrm{d}}y+Ce^{-2\gamma s}.
\end{eqnarray}
Considering  $J(w(s),s)$ defined in \eqref{f3}, we obtain from
equation (\ref{eqw}) and integration by parts
\begin{multline}\label{t1}
e^{\gamma s}\frac{d}{ds}J(w(s),s)= - \iint (\ps
w)^2\rho{\mathrm{d}}y + \iint(\py w)^2(1-y^2)\rho{\mathrm{d}}y
+\frac{2p+2}{(p-1)^2}\iint w^2\rho{\mathrm{d}}y\nonumber\\
 - \iint |w|^{p+1}\rho{\mathrm{d}}y+
(\gamma +\frac{p+3}{p-1}-2N) \iint w \ps w (s) \rho{\mathrm{d}}y-2\iint w\ps wy \rho'{\mathrm{d}}y\qquad\qquad\\
-2\iint \ps w\py w
y\rho{\mathrm{d}}y-e^{-\frac{2ps}{p-1}}\iint
wf\Big(e^{\frac{2s}{p-1}}w\Big){\rho}{\mathrm{d}}y-(N-1) e^{-s} \iint  w\py w \frac{\rho}r {\mathrm{d}}y\nonumber\\
 - e^{-\frac{2ps}{p-1}}\iint w
g\Big(r_0+ye^{-s},T_0-e^{-s},e^{\frac{(p+1)s}{p-1}}\partial_yw,e^{\frac{(p+1)s}{p-1}}(\partial_sw+y\partial_y
w+\frac{2}{p-1}w)\Big){\rho}{\mathrm{d}}y.
\end{multline}
Combining (\ref{defF}), (\ref{f2}) and (\ref{t1}),   we write
\begin{eqnarray}\label{theta}
&e^{\gamma s}\frac{d}{ds}J(w(s),s) \le \frac{p+3}2
\big(E_0(w(s))+I(w(s),s)\big) -\frac{p-1}{4} \iint (\py
w)^2(1-y^2)\rho{\mathrm{d}}y&\nonumber\\
 &-\frac{p+1}{2(p-1)}\iint w^2\rho{\mathrm{d}}y
-\frac{p-1}{2(p+1)}\iint |w|^{p+1}\rho{\mathrm{d}}y&\nonumber\\
&+
\underbrace{(\gamma +\frac{p+3}{p-1}-2N+\frac{p+3}2 e^{-\gamma s}) \iint w\ps w \rho{\mathrm{d}}y}_{J_1(s)}&\\
&\underbrace{+\frac8{p-1} \iint w\ps w\frac{y^2}{1-y^2}
\rho{\mathrm{d}}y}_{J_2(s)} \underbrace{-2\iint \ps w\py w
y\rho{\mathrm{d}}y}_{J_3(s)}
 \underbrace{-e^{-\frac{2ps}{p-1}}\iint
wf\Big(e^{\frac{2s}{p-1}}w
\Big){\rho}{\mathrm{d}}y}_{J_4(s)}&\nonumber\\ &\underbrace{  -
e^{-\frac{2ps}{p-1}}\iint w
g\Big(r_0+ye^{-s},T_0-e^{-s},e^{\frac{(p+1)s}{p-1}}\partial_yw,e^{\frac{(p+1)s}{p-1}}(\partial_sw+y\partial_y
w+\frac{2}{p-1}w)\Big) {\rho}{\mathrm{d}}y}_{J_5(s)}&
\nonumber\\
&+ \underbrace{\frac{p+3}2
e^{-\frac{2(p+1)s}{p-1}}\displaystyle\iint
F(e^{\frac{2}{p-1}s}w)\rho {\mathrm{d}}y}_{J_6(s)}\underbrace{-(N-1)
e^{-s} \iint w\py w \frac{\rho}r
{\mathrm{d}}y}_{J_7(s)}.\nonumber
\end{eqnarray}
We now study each of the last five terms. To estimate $J_1(s)$ and
$J_2(s)$, we use   the Cauchy-Schwartz inequality to have
\begin{eqnarray}\label{J1}
 |J_1(s)|& \le&
 Ce^{\frac{\gamma s}2}\iint (\ps w)^2\frac{\rho}{1-y^2}{\mathrm{d}}y+
C e^{-\frac{\gamma s}2}\iint w^2\rho
{\mathrm{d}}y.\qquad\qquad
\end{eqnarray}
\begin{eqnarray*}\label{J2}
|J_2(s)| &\le&
 Ce^{\frac{\gamma s}2}\iint (\ps w)^2\frac{\rho}{1-y^2}{\mathrm{d}}y+
C  e^{-\frac{\gamma s}2}\iint w^2\frac{y^2\rho}{1-y^2}
{\mathrm{d}}y.
\end{eqnarray*}
Recalling the following Hardy-Sobolev estimate (see Appendix B
page 1163 in \cite{MZajm03} for the proof):
\begin{equation}\label{hs}
\iint h^2 \frac \rho{1-y^2} {\mathrm{d}}y \le C\iint h^2 \rho
{\mathrm{d}}y + C \iint (h'(y))^2 \rho(1-y^2){\mathrm{d}}y,
\end{equation}
we  conclude that
\begin{eqnarray}\label{J22}
|J_2(s)| &\le& C e^{\frac{\gamma s}2} \iint (\ps
w)^2\frac{\rho}{1-y^2}{\mathrm{d}}y+ C e^{-\frac{\gamma
s}2}\iint w^2\rho {\mathrm{d}}y\nonumber\\
&&+ C  e^{-\frac{\gamma
s}2}\iint \!(\py w)^2\rho(1-y^2) {\mathrm{d}}y.
\end{eqnarray}
Using the Cauchy-Schwartz inequality, we have
\begin{eqnarray}\label{J3}
|J_3(s)| \le
 C e^{\frac{\gamma s}2}\iint (\ps w)^2\frac{\rho}{1-y^2}{\mathrm{d}}y+ C  e^{-\frac{\gamma s}2}\iint \!(\py w)^2\rho(1-y^2) {\mathrm{d}}y.
\end{eqnarray}
From \eqref{born}, we write
\begin{eqnarray}\label{J4}
|J_4(s)|+|J_6(s)| &\le&C e^{-\gamma s}+
 C e^{-\gamma s}\iint |w|^{p+1}\rho{\mathrm{d}}y.
\end{eqnarray}
In a similar way,    using the hypothesis  $(H_g) $ and (\ref{hs}),
we have
\begin{eqnarray}\label{J5}
|J_5(s)|&\le& Ce^{-\gamma s} \iint (\ps w)^2\frac{\rho}{1-y^2}
{\mathrm{d}}y +C e^{-\gamma s}\iint
(\py w)^2\rho(1-y^2){\mathrm{d}}y\nonumber\\
&& +Ce^{-\gamma s}\iint w^2\rho {\mathrm{d}}y+ C e^{-\gamma
s}.
\end{eqnarray}
Using \eqref{range} and  (\ref{hs}), we write
\begin{eqnarray}\label{CS}
|J_7(s)| \le Ce^{-\frac{s}2}\iint (\py w)^2\rho(1-y^2)
{\mathrm{d}}y +Ce^{-\frac{s}2} \iint w^2 {\rho} {\mathrm{d}}y.
\end{eqnarray}
Finally,  using (\ref{theta}), (\ref{J1}), (\ref{J22}), (\ref{J3}),
(\ref{J5}) and (\ref{CS}), we deduce that
\begin{eqnarray}\label{Jprime}
e^{\gamma s}\frac{d}{ds}J(w(s),s) &\le &\frac{p+3}2
\Big(E_0(w(s))+I(w(s),s)\Big)\\
&&+\Big(C e^{-\frac{\gamma s}2}-\frac{p-1}{4}\Big) \iint (\py w)^2(1-y^2)\rho{\mathrm{d}}y\nonumber\\
&&+\Big(C e^{-\frac{\gamma s}2}-\frac{p+1}{2(p-1)}\Big)\iint
w^2\rho{\mathrm{d}}y\nonumber\\
&&+\Big(C e^{-\frac{\gamma s}2}-\frac{p-1}{2(p+1)}\Big)\iint |w|^{p+1}\rho{\mathrm{d}}y\nonumber\\
&&+C e^{\frac{\gamma s}2} \iint (\ps
w)^2\frac{\rho}{1-y^2}{\mathrm{d}}y+C e^{-\gamma s}.\nonumber
\end{eqnarray}
 From  (\ref{E}) and (\ref{Jprime}), we obtain
\begin{eqnarray}
 \frac d{ds}E(w(s),s)&\le &Ce^{-2\gamma s}+\frac{p+3}2 e^{-\gamma s}E(w(s),s)\nonumber\\
 &&+\Big(Ce^{-\frac{\gamma
s}2}-\frac 4{p-1}\Big) \iint (\ps w)^2 \frac{\rho}{1-y^2} {\mathrm{d}}y \nonumber\\
&&+\Big(C e^{-\frac{\gamma s}2}-\frac{p+1}{2(p-1)}\Big)e^{-\gamma
s}\iint w^2\rho{\mathrm{d}}y
\nonumber\\
&& +\Big(C e^{-\frac{\gamma s}2}-\frac{p-1}{4}\Big) e^{-\gamma s}\iint  (\py w)^2(1-|y|^2) \rho {\mathrm{d}}y\\
&& +\Big(C e^{-\frac{\gamma s}2}-\frac{p-1}{2(p+1)}\Big)e^{-
\gamma s}\iint |w|^{p+1}\rho{\mathrm{d}}y . \nonumber
\end{eqnarray}
We now choose $S_0\ge 0$, large enough, so that for all $s\ge
S_0$, we have
\begin{eqnarray*}
 \frac{p-1}{4}-C e^{-\frac{\gamma s}2} \ge 0, \
 \frac{p+1}{2(p-1)}
-C e^{-\frac{\gamma s}2} \ge 0, \ \frac{p-1}{2(p+1)}-C
e^{-\frac{\gamma s}2}\ge 0,\  \frac 1{p-1}-\C e^{-\frac{\gamma
s}2}\ge0.
\end{eqnarray*}
Then, we deduce that, for all $s\ge \max(S_0,-\log T(r_0),-\log
\frac{r_0}2, -4\log r_0)$, we have
\begin{eqnarray}\label{E111}
 \frac d{ds}E(w(s),s)\le Ce^{-2\gamma s}+\frac{p+3}2 e^{-\gamma s}E(w(s),s)-\frac 3{p-1} \iint (\ps w)^2 \frac{\rho}{1-y^2}
 {\mathrm{d}}y.\qquad
\end{eqnarray}
 This yields (i) of Proposition
\ref{prophamza}.

\medskip

(ii)
\no We finish  the proof  of  Proposition \ref{prophamza} here. More
precisely, we
prove that there exists $S_1(p,N,M,q)\in \R$ such that,  for all $ x_0\in \R^N$ and $ T_0\in (0,T(x_0)]$,
\begin{equation}\label{254}
\forall \ s\ge  \max\left(-\log T_0, S_1,\right)
\ \ \ H(w_{x_0,T_0}(s),s)\ge 0.
\end{equation}
 We  give the proof only  in the case where $x_0$ is a non
characteristic point. Note that the case where $x_0$ is a
 characteristic point  can be  done  exactly as in Appendix A page 119 in \cite{MZjfa07}.
 If $x_0$ is a non
characteristic point,  the argument is the
same as in the corresponding part in \cite{AMimrn01}. We write the
proof for completeness. Arguing by contradiction, we assume that
there exists  a non characteristic point $ x_0\in \R^N$, $ T_0\in
(0,T(x_0)]$  and  $ s_1\ge \max\left(-\log T_0, S_1,\right)$ such that $H(w(s_1),s_1)<0$, where
$w=w_{x_0,T_0}$.
 By
definition \eqref{defH} of $H$, we write
\begin{eqnarray*}
&&H(W(s),s)\ge\mu e^{-2\gamma s} -{e^{\frac{p+3}{2}e^{-\gamma s}}}\Big(\frac1{p+1}+Ce^{-2\gamma s}\Big)\iint |W|^{p+1} \rho {\mathrm{d}}y\\
&& +e^{\frac{p+3}2e^{-\gamma s}}\left(\left(\frac 12 -C e^{-\gamma
s}\right)\iint (\ps W)^2 \rho {\mathrm{d}}y
+\left(\frac{p+1}{(p-1)^2} - Ce^{-\gamma s}\right)\iint W^2 \rho {\mathrm{d}}y\right)\\
&\ge & -\frac 2{p+1}\iint |W|^{p+1} \rho {\mathrm{d}}y.
\end{eqnarray*}
if $s\ge S_2(p,N,q,M)\ge S_1(p,N,q,M)$ for some $S_2(p,N,q,M)\in\R$ large enough.
Using this inequality together with the fact that $H(W(s),s)$ is
decreasing by the remark following Proposition \ref{prophamza}, we
see that the argument used by Antonini and Merle in Theorem 2 page
1147 in \cite{AMimrn01} for the equation \eqref{par} works here and
we get the blow-up criterion. This concludes the proof of
Proposition \ref{prophamza}.\Box

\section{Blow-up results related to non-characteristic points}\label{secnonchar}

Let us first introduce for all $|d|<1$ the following solitons defined by
\begin{equation}\label{defkd}
\kappa(d,y)=\kappa_0 \frac{(1-d^2)^{\frac 1{p-1}}}{(1+dy)^{\frac 2{p-1}}}\mbox{ where }\kappa_0 =
\left(\frac{2(p+1)}{(p-1)^2}\right)^{\frac 1{p-1}} \mbox{ and }|y|<1.
\end{equation}
Note that $\kappa(d)$ is a stationary solution of \eqref{eqw}, in
the particular case where $(f,g)\equiv (0,0)$ and in one space
dimension.

\medskip

\noindent  Adapting the analysis of \cite{MZjfa07} and \cite{MZcmp08}, we claim the following:
\begin{theor}[Blow-up behavior and regularity of the blow-up set on $\RR$]\label{thbb}$ $\\
(i) {\bf (Regularity related to $\RR$)} $\RR\neq \emptyset$, $\RR\cap \R^*_+$ is an open set, and $x\mapsto T(x)$ is of class $C^1$ on $\RR\cap \R^*_+$.\\
(ii) {\bf (Blow-up behavior in similarity variables)} There exist $\mu_0>0$ and $C_0>0$ such that for all $r_0\in\RR\cap \R^*_+$, there exist
$\theta(r_0)=\pm 1$ and $s_0(r_0)\ge - \log T(r_0)$ such that for all $s\ge s_0$:
\begin{equation}\label{profile}
\left\|\vc{w_{r_0}(s)}{\partial_s w_{r_0}(s)}-\theta(r_0)\vc{\kappa(T'(r_0))}{0}\right\|_{\H}\le C_0 e^{-\mu_0(s-s_0)}.
\end{equation}
Moreover, $E_0(w_{r_0}(s)) \to E_0(\kappa_0)$ as $s\to \infty$.
\end{theor}
{\bf Remark}: As stated in the introduction, this result holds also when $N=1$, with  no symmetry assumptions an  initial data, for all $r_0\in \R$, even $r_0=0$; when $N\ge 2$ and if $0\in\RR$, the asymptotic behavior of $w_0$ remains open.\\
{\it Proof}: As in the non-perturbed radial case (take
$(f,g)\equiv(0,0)$ in \eqref{equ}) treated in \cite{MZbsm11}, we
need to make some minor adaptations to  the one-dimensional
non-perturbed case treated in \cite{MZjfa07} and \cite{MZcmp08}. It
happens that the same adaptation pattern  works in the present case,
and that is the reason why we don't mention it, and refer the
reader to the proof of Theorem 1 in page 358 of \cite{MZbsm11}. The
only points to check are the following:

- {\it Continuity with respect to the scaling parameter}: Due to the fact that equation \eqref{equ} is no longer invariant
under the scaling
\[
\lambda \mapsto u_\lambda(\xi, \tau) = \lambda^{\frac 2{p-1}}u(\lambda \xi, \lambda \tau),
\]
we need to understand the continuous dependence of the solutions of the following family of equations
\begin{eqnarray}\label{eqxl}
\partial^2_{t} u =\partial^2_r u+\frac{(N-1)\lambda}{x+\lambda r} \pr u+|u|^{p-1}u
+\lambda^{\frac{2p}{p-1}}f\Big(\lambda^{\frac{-2}{p-1}}u\Big)\nonumber\\ +
\lambda^{\frac{2p}{p-1}}g\Big(\lambda r,\lambda
t,\lambda^{-\frac{p+1}{p-1}}\partial_ru,\lambda^{-\frac{p+1}{p-1}}\partial_tu\Big),
\end{eqnarray}
with respect to initial data and the parameters $x\ge 0$ and $\lambda>0$ (including the limit as $\lambda \to 0$ 
 and this is a classical estimate.

\medskip

- {\it A new statement for the trapping result}:
This
is due to the fact that equation \eqref{eqw} in similarity variables depends on
a parameter $r_0>0$ and contains new terms of order $e^{-\gamma s}$ ($\gamma$ is defined in \eqref{gamma}) (it is no longer autonomous).
This is the trapping result in
our setting:
\begin{theor}{\bf (Trapping near the set of non zero stationary solutions of \eqref{eqw})}\label{thtrap}
\\
For all $\rho_0>0$, there exist positive $\epsilon_0$,
$\mu_0$ and $C_0$ such that for all $\epsilon^*\le \epsilon_0$, there exists
$s_0(\epsilon^*)$ such that if $r_0\ge \rho_0$, $s^*\ge s_0$ and $w\in C([s^*, \infty), \H)$ is a solution of equation \eqref{eqw} with
\begin{equation}\label{highenergy}
\forall s\ge s^*,\;\;E(w(s),s)\ge E_0(\kappa_0)-e^{-\frac {\gamma s}2},
\end{equation}
and
\begin{equation*}
\left\|\vc{w(s^*)}{\partial_s w(s^*)}-\omega^*\vc{\kappa(d^*,\cdot)}{0}\right\|_{\H}\le \epsilon^*
\end{equation*}
for some $d^*\in(-1,1)$ and $\omega^*=\pm 1$,
then there exists $d_\infty\in (-1, 1)$ such that
\[\left|\argth{d_\infty} - \argth{d^*}\right|\le C_0 \epsilon^*,\]
and for all $s\ge s^*$,
\begin{equation*}
\left\|\vc{w(s)}{\partial_s w(s)}-\omega^*\vc{\kappa(d_\infty, \cdot)}{0}\right\|_{\H}\le C_0 \epsilon^* e^{-\mu_0(s-s^*)}.
\end{equation*}
\end{theor}
{\it Proof}: The proof follows the pattern of the radial case treated in \cite{MZbsm11}. For that reason, we refer the reader to the Proof of Theorem 2 page 360 in that paper, and focus in the following only on how to treat the new terms generated by the perturbations $f$ and $g$ in \eqref{gen}. With respect to the pure power case in one space dimension, the difference comes from the linearization of \eqref{eqw} around the stationary solutions $\kappa(d,y)$ in \eqref{defkd}, where we see the following lower order terms:
\begin{eqnarray}
&&\left|\frac{(N-1)e^{-s}}{r_0+ye^{-s}}\py w\right|\le \frac
2{\rho_0} (N-1) e^{-s}|\py w|;\label{small}\\
&&e^{-\frac{2ps}{p-1}}\left|f\left(e^{\frac{2s}{p-1}}w\right)\right|
\le CMe^{-\frac{2(p-q)s}{p-1}}
+CMe^{-\frac{2(p-q)s}{p-1}}\big|w\big|^p;\nonumber\\
&&e^{-\frac{2ps}{p-1}}
\Big|g\Big(r_0+ye^{-s},T_0-e^{-s},e^{\frac{(p+1)s}{p-1}}\partial_yw,e^{\frac{(p+1)s}{p-1}}(\partial_sw+y\partial_y
w+\frac{2}{p-1}w)\Big)\Big|\nonumber\\
&&\le CMe^{-s}\Big(1+\big|\partial_sw\big|+\big|\partial_y
w\big|+\big|w\big|\Big).\nonumber
\end{eqnarray}
as soon as $r_0\ge \rho_0>0$, and $s\ge -\log \frac{\rho_0}2$.
For more details an the adaptation, we refer the reader to the proof
of  Theorem 1 in page 358 of \cite{MZbsm11}.\Box

\section{Blow-up results related to characteristic points}\label{secchar}
The first question in this case is of course the existence of examples of initial data with $\SS \not= \emptyset$.
If the perturbation $g$ introduced in \eqref{gen} does not depend on $|x|$, then the existence of such an example
 follows from the knowledge of the blow-up behavior at non-characteristic points, as in the pure power nonlinearity
 case \eqref{par}. If $g$ depends on $|x|$, then we need to apply the constructive method of C\^ote and Zaag \cite{CZarxiv},
  which relies fundamentally on the knowledge of the blow-up behavior near a characteristic point. For that reason, we leave
  the existence issues to the end of the section, and start with the description of the blow-up features near characteristic
   points. More precisely, we proceed in two sections:\\
- In Section \ref{description}, we consider arbitrary blow-up solutions having a non-zero characteristic point, and we give
a full description of its blow-up behavior and its blow-up set near this characteristic point.\\
- In Section \ref{existence}, we prove the existence of such a solution, and also give some criteria for the existence
or the non-existence of characteristic points.
\subsection{Description of the blow-up behavior and the blow-up set near a characteristic point}\label{description}

Now, given $r_0\in\SS \cap \R^*_+$, we have the same description for
the asymptotic  of $w_{r_0}$ as in the one-dimensional case with no
perturbations (i.e. for equation (\ref{equ}) with $(f,g)\equiv(0,0)$)
refined recently by C\^ote and Zaag  in  \cite{CZarxiv}. In order to state the result, let us introduce
\begin{equation}\label{solpart}
\bar \zeta_i(s) = \left(i-\frac{(k+1)}2\right)\frac{(p-1)}2\log s + \bar\alpha_i(p,k)
\end{equation}
where the sequence $(\bar\alpha_i)_{i=1,\dots,k}$ is uniquely determined by the fact that $(\bar \zeta_i(s))_{i=1,\dots,k}$ is an explicit solution with zero center of mass for the following ODE system:
\begin{equation} \label{eq:tl}
 \frac 1{c_1}\dot \zeta_i = e^{ - \frac{2}{p-1} (\zeta_i - \zeta_{i-1}) } - e^{- \frac{2}{p-1} (\zeta_{i+1} - \zeta_i) },
\end{equation}
where $c_1=c_1(p)>0$ and $\zeta_0(s)\equiv \zeta_{k+1}(s) \equiv 0$ (see Section 2 in \cite{CZarxiv} for a proof of this fact). Note that $c_1=c_1(p)>0$ is a constant appearing in system \eqref{eqz}, itself
inherited from Proposition 3.2 of \cite{MZajm10}.
With this definition, we can state our result
 (for the statement in one space dimension, see Theorem 6 in
\cite{MZajm10} and  Theorem 1 in \cite{CZarxiv}):
\begin{theor}\label{new}
{\bf (Description of the behavior of $w_{r_0}$ where $r_0$ is
characteristic)} Consider
$r_0\in \SS\cap \R^*_+$. Then, there is $\zeta_0(r_0)\in \R$
such that
\begin{equation}\label{cprofile00}
\left\|\vc{w_{r_0}(s)}{\ps w_{r_0}(s)}
 - \theta_1\vc{\d\sum_{i=1}^{k(r_0)} (-1)^{i+1}\kappa(d_i(s),\cdot)}0\right\|_{\H} \to 0\mbox{ and }E_0(w_{r_0}(s))\to k(r_0)E_0(\kappa_0)
\end{equation}
as $s\to \infty$, for some
\begin{eqnarray}
&&k(r_0)\ge 2,\label{pb}\\
&&\theta_i=\theta_1(-1)^{i+1},\ \theta_1=\pm 1\label{ei}
\end{eqnarray}
 and continuous $d_i(s) = -\tanh \zeta_i(s)$ with
\begin{equation}\label{equid}
\zeta_i(s) = \bar \zeta_i(s) + \zeta_0,
\end{equation}
where $\bar \zeta_i(s)$ is introduced above in \eqref{solpart}.
\end{theor}
{\bf Remark}:
As stated in the introduction, this result holds also when $N=1$, with  no symmetry assumptions an  initial data, for all $r_0\in \R$, even $r_0=0$; when $N\ge2$ and if $0\in\SS$, the asymptotic behavior of $w_0$ remains open.\\
{\it Proof}: As in the one-dimensional case with no perturbations
(i.e. for equation (\ref{equ}) with $(f,g)\equiv(0,0)$),
 the proof of the asymptotic behavior and the geometric
results on $\SS$ (see Theorem \ref{thgeo} below) go side by side.
Evidently the refined  description given by \eqref{equid} is
obtained as in \cite{CZarxiv}. We leave the proof after the
statement of Theorem \ref{thgeo}.

\Box

Let us note that
we get the following result on the energy behavior from the asymptotic behavior
at a non-characteristic point (see (ii) of Theorem \ref{thbb}) and at a characteristic point (see Theorem \ref{new}):
\begin{coro}[A criterion for non-characteristic points]\label{corcriterion}$ $\\
 For all $r_0>0$, there exist $C_3(r_0)>0$ and $S_3(r_0)\in\R$ such that:\\
(i) For all $r\in[\frac{r_0}2, \frac{3r_0}2]$ and $s\ge S_3$, we have
\[
E_0(w_r(s))\ge k(r)E_0(\kappa_0)-C_3(r_0)e^{-\gamma s}.
\]
(ii) If for some $r\in[\frac{r_0}2, \frac{3r_0}2]$ and $s\ge S_3$, we have
\[
E_0(w_r(s))<2 E_0(\kappa_0)-C_3(r_0)e^{-\gamma s},
\]
then $r\in \RR$.
\end{coro}
{\bf Remark}: With respect to the statement in  one-space dimensions with no perturbations,  (Corollary 7 in \cite{MZajm10}), this statement has additional exponentially small terms. This comes from the fact that the functional $E(w(s),s)$ is no longer decreasing, and that one has to work instead with the functional $H(w(s),s)$ \eqref{defH} which is decreasing, and differs from $E(w(s),s)$ by exponentially small terms, uniformly controlled for $r\in[\frac{r_0}2, \frac{3r_0}2]$ thanks to the uniform estimates of Proposition \ref{boundedness}.\\
{\it Proof}: If one replaces $E(w(s),s)$ by $H(w(s),s)$, then the proof is
straightforward from Theorems \ref{thbb} and \ref{new} together with
the monotonicity of $H(w(s),s)$ (see \eqref{defH} and
\eqref{prophamza}). Since the difference between the two functionals
is exponentially small, uniformly for $r\in[\frac{r_0}2,
\frac{3r_0}2]$ (see \eqref{defH}, \eqref{edoF} and Proposition
\ref{boundedness}), we get the conclusion of Corollary
\ref{corcriterion}.
\Box

\bigskip

Finally, we give in the following some geometric information related to
characteristic points (for the statement in one space dimension, see Theorem 1,
Theorem 2 and the following remark in \cite{MZisol10}):
\begin{theor}\label{thgeo}{\bf (Geometric considerations on $\SS$)}\\ 
(i) {\bf (Isolatedness of characteristic points)}  Any $r_0\in \SS\cap \R^*_+$  is isolated.\\
(ii) {\bf (Corner shape of the blow-up curve at characteristic
points)} If $r_0\in \SS\cap \R^*_+$  with  $k(r_0)$ solitons
and $\zeta_0(r_0)\in \R$ as center of mass of the solitons' center
as shown in \eqref{cprofile00} and \eqref{equid}, then
\begin{eqnarray}
T'(r)+\theta(r)&\sim&\frac{\theta(r)\nu e^{-2\theta(r)\zeta_0(r_0)}}{|\log|r-r_0||^{\frac{(k(r_0)-1)(p-1)}2}}\label{cor1}\\
T(r)-T(r_0)+|r-r_0|&\sim&\frac{\nu
e^{-2\theta(r)\zeta_0(r_0)}|r-r_0|}{|\log|r-r_0||^{\frac{(k(r_0)-1)(p-1)}2}}\label{cor0}
\end{eqnarray}
as $r\to r_0$, where $\theta(r) = \frac{r-r_0}{|r-r_0|}$ and
$\nu=\nu(p)>0$.
\end{theor}
{\it Proof}: See below.\\
{\bf Remark}:
As stated in the remark after Theorem 3, our result holds for $N=1$, with  no symmetry assumptions an  initial data, for all $r_0\in \R$, even $r_0=0$; when $N\ge2$, and  if $0\in\SS$, the asymptotic behavior of $w_0$ remains open.\\
{\it Proof}: As in the one-dimensional case with no perturbations

{\bf Remark}: Note from (i) that the multi-dimensional version $U(x,t)=u(|x|,t)$ has a finite number
 of concentric spheres of characteristic points in the set $\{\frac 1R < |x|<R\}$ for every $R>1$.
 This is consistent with our conjecture in \cite{MZisol10} where we guessed that in dimension $N\ge 2$,
 the $(N-1)$-dimensional Hausdorff measure of $\SS$ is bounded in compact sets of $\R^N$. Note that this
 conjecture is related to the result of Vel\'azquez who proved in \cite{Viumj93}
 that the $(N-1)$-dimensional Hausdorff measure of the blow-up set for the semilinear heat
  equation with subcritical power nonlinearity is bounded in compact sets of $\R^N$.

\bigskip

As a consequence of our analysis, particularly the lower bound on $T(r)$ in \aref{cor0},
 we have the following estimate on the blow-up speed in the backward light cone with
 vertex $(r_0, T(r_0))$ where $r_0> 0$ (for the statement in one space dimension, see Corollary 3 in \cite{MZisol10}):
\begin{coro}\label{corspeed}{\bf (Blow-up speed in the backward light cone)} For all $r_0>0$,
there exists $C_4(r_0)>0$ such that for all $t\in[0, T(r_0))$, we have
\[
\frac{|\log(T(r_0)-t)|^{\frac{k(r_0)-1}2}}{C_4(r_0)(T(r_0)-t)^{\frac 2{p-1}}}\le
\sup_{|x-r_0|<T(r_0)-t}|u(x,t)|\le \frac{C_4(r_0) |\log(T(r_0)-t)|^{\frac{k(r_0)-1}2}}{(T(r_0)-t)^{\frac 2{p-1}}}.
\]
\end{coro}
{\bf Remark}: Note that when $r_0\in\RR\cap \R^*_+$, the blow-up rate of $u$ in the backward light cone with vertex $(r_0, T(r_0))$ is given by the solution of the associated ODE $u"=u^p$. When $r_0\in\SS\cap \R^*_+$, the blow-up rate is higher and quantified, according to $k(r_0)$, the number of solitons appearing in the decomposition \eqref{cprofile00}.\\
{\it Proof}: When $r_0\in\RR$, the result follows from the fact that the convergence in \aref{profile} is true also in $L^\infty\times L^2$ from
\eqref{profile}
and the Sobolev embedding in one dimension. When $r_0\in\SS$, see the proof of Corollary 3 of \cite{MZisol10} given in Section 3.3 of that paper.

\Box

\bigskip

{\it Proof of Theorems \ref{new} and \ref{thgeo}}: The proof follows the pattern of the original proof, given in \cite{MZjfa07}, \cite{MZajm10}, \cite{MZisol10} and \cite{CZarxiv}. In the following, we recall its different parts.

\bigskip

 {\bf Part 1: Proof of \eqref{cprofile00} without \eqref{pb} nor \eqref{ei} and with the estimate
\begin{equation}\label{decuple}
\zeta_{i+1}(s)-\zeta_i(s) \to \infty\mbox{ as }s\to\infty
\end{equation}
 instead of \eqref{equid}} (note that
\eqref{decuple}
is
meaningful only when $k(r_0)\ge 2$).

 The original statement of this part is given in Theorem 2 (B) page 47 in \cite{MZjfa07} and the proof in section 3.2 page 66 in that paper. Note that this part doesn't exclude the possibility of having $k(r_0)=0$ or $k(r_0)=1$. The adaptation is straightforward. As in the non-characteristic case above, one has to use the Duhamel formulation in the radial which may be derived from \cite{SSnyu98}.


\bigskip

\label{pagepart2}{\bf Part 2: Assuming that \eqref{pb} is true, we prove \eqref{ei} with the estimates
\begin{eqnarray}
|\zeta_i(s)-\bar \zeta_i(s)|&\le & C,\label{equidnew}\\
T(r)-T(r_0)+|r-r_0|&\le &\frac{C|r-r_0|}{|\log|r-r_0||^{\frac{(k(r_0)-1)(p-1)}2}}\nonumber
\end{eqnarray}
instead of \eqref{equid} and \eqref{cor0}}.

The original statement is given in Propositions 3.1 and 3.13 in \cite{MZajm10}. The reader has to read Section 3 and Appendices B and C in that paper.
The adaptation is straightforward, except for the effect of the new terms in equation \eqref{eqw},
which produce exponentially small terms in many parts of the proof (see \eqref{small}). In particular, Lemma 3.11 of \cite{MZajm10} has to be changed by adding $Ce^{-\gamma s}$ where   $\gamma$  is defined in \eqref{gamma} to the right of all the differential inequalities.

\bigskip

 {\bf Part 3: Proof of \eqref{pb} and the fact that the interior of $\SS$ is empty}.

The original statement is given in Proposition 4.1 of \cite{MZajm10}. The adaptation is as delicate as in \cite{MZbsm11}.
In particular, it involves the ruling-out of the occurrence of the case where, locally near the origin, the blow-up set of the multi-dimensional version $U(x,t)$ is a forward light cone with vertex $(0,T(0))$
(see Lemma 4.5 page 367 in \cite{MZbsm11}). As in that paper, the proof of the non-occurrence of this case is based in particular on a local energy estimate by Shatah and Struwe \cite{SSnyu98}. For the reader's convenience, we adapt in Appendix \ref{appenergy} that energy estimate to our case \eqref{gen}, namely to perturbations of the pure power equation \eqref{par}. For the other arguments, we refer to  the corresponding part in \cite{MZbsm11} (see Part 3 page 336 in that paper).

\bigskip

{\bf Part 4: Proof of Theorem \ref{thgeo} with \eqref{cor1} and \eqref{cor0} replaced by
\begin{eqnarray*}
\d\frac{1}{C_0|\log(r-r_0)|^{\frac{(k(r_0)-1)(p-1)}2}}\le& T'(r)+\frac{r-r_0}{|r-r_0|} &\le \frac{C_0}{|\log(r-r_0)|^{\frac{(k(r_0)-1)(p-1)}2}},\\
\d\frac{|r-r_0|}{C_0|\log(r-r_0)|^{\frac{(k(r_0)-1)(p-1)}2}}\le &T(r)- T(r_0)+|r-r_0|& \le \frac{C_0|r-r_0|}{|\log(r-r_0)|^{\frac{(k(r_0)-1)(p-1)}2}}.
\end{eqnarray*}
}

The analogous statement in one space dimension with no perturbations is given in Theorems 1 and 2
in \cite{MZisol10}. Thus, one needs to say how to adapt the analysis of the paper \cite{MZisol10} to the present case. As in \cite{MZbsm11}, three ingredients are needed in the proof:\\
- the trapping result stated in Theorem \ref{thtrap};\\
- the energy criterion stated in Corollary \ref{corcriterion};\\
- the dynamics of equation \eqref{eqw} around a decoupled sum of solitons performed in \cite{MZisol10} and presented in Part 3 above.
Note that we have already adapted all these ingredients to the present context. With this fact, the adaptation given in \cite{MZbsm11} works here. See Part 4 page 371 in that paper for more details.

\bigskip

{\bf Part 5: Proof of \eqref{equid}, \eqref{cor1} and \eqref{cor0}}

This part corresponds to the contributions brought in \cite{CZarxiv} in the one-dimensional case. The orginal statements in the one-dimensional case are given in Theorem 1.1 and Corollary 1.4 in that paper. Following Part 2 where we proved that \eqref{cprofile00} holds with \eqref{equid} replaced by \eqref{equidnew}, a crucial step in one-space dimension was to prove that the solitons' centers satisfy the following ODE system for $s$ large enough:
\begin{equation}\label{eqz}
 \frac 1{c_1}\dot \zeta_i = e^{ - \frac{2}{p-1} (\zeta_i - \zeta_{i-1}) } - e^{- \frac{2}{p-1} (\zeta_{i+1} - \zeta_i) }+O\left(\frac 1{s^{1+\eta}}\right)
\end{equation}
for some $\eta>0$. In \cite{CZarxiv}, we were able to use some ODE tools (particularly the Lyapunov convergence theorem) to further refine estimate \eqref{equidnew} and prove that
\[
\zeta_i(s) = \bar \zeta_i(s) + \zeta_0+o\left(\frac 1{s^\eta}\right)\mbox{ as }s \to \infty.
\]
Since we have for all $|d_1|<1$ and $|d_2|<1$
\[
\|\kappa(d_1)-\kappa(d_2)\|_{\H}\le C|\argth d_1 - \argth d_2|
\]
(see estimate (174) page 101 in \cite{MZjfa07} for a proof of this fact),
estimate \eqref{cprofile00} remains unchanged if one slightly modifies $\zeta_i(s)$ by setting $\zeta_i(s) = \bar \zeta_i(s) + \zeta_0$ which is the desired estimate in \eqref{equid}. That was the argument in one space dimension.\\
 In our setting, since our perturbative terms contribute with additional exponentially decaying terms to the equation (see \eqref{small} and Part 2 above), we obtain that $\zeta_i(s)$ satisfy the same ODE system \eqref{eqz}. Thus, the refinements of \cite{CZarxiv} hold here with no need for any further adaptations, and \eqref{equid} holds.

\medskip

As for estimates \eqref{cor1} and \eqref{cor0}, let us point out that in one space dimension, they are derived in \cite{CZarxiv} as direct consequences of \eqref{equid} on the one hand, and on the other hand a small improvement of the last argument of the paper \cite{MZisol10} based on the equation in similarity variables. Since in our setting, \eqref{equid} holds and the equation in similarity variables differs from the one dimensional case with exponentially decaying terms (see \eqref{small}), the same argument holds. See Section 2 in \cite{CZarxiv}.
\subsection{Existence and non-existence of characteristic points}\label{existence}
Proceeding as in \cite{CZarxiv}, we have the following result:
\begin{theor}\label{mainth} {\bf (Existence of a solution with prescribed blow-up behavior at a characteristic point)} For any $r_0>0$ and $k\ge 2$,
there exists a blow-up solution $u(r,t)$ to equation \eqref{equ}
with $r_0\in\SS$ such that
\begin{equation}\label{cprofile0}
\left\|\vc{w_{r_0}(s)}{\ps w_{r_0}(s)} - \vc{\ds\sum_{i=1}^{k} (-1)^{i+1}\kappa(d_i(s))}0\right\|_{\H} \to 0\mbox{ as }s\to \infty,
\end{equation}
with
\begin{equation}\label{refequid1}
d_i(s) = -\tanh \zeta_i(s), \quad
\zeta_i(s) = \bar \zeta_i(s) + \zeta_0
\end{equation}
for some $\zeta_0\in \R$, where $\bar \zeta_i(s)$ is defined in \eqref{solpart}.
\end{theor}
{\bf Remark}: When $N=1$, we can take $r_0=0$. When $N\ge 2$, the multi-dimensional version $U(x,t)=u(|x|,t)$ has a sphere of characteristic points. Note also that this result uses the same argument as for Theorem \ref{new}, in particular, the analysis of the ODE system  \eqref{eqz}.
If we simply want an argument for the existence of a blow-up solution with a characteristic point without caring about the number of solitons, then we have a more elementary proof which holds, however, only when $g$ does not depend on $|x|$. See the remark following Theorem \ref{thexis} below.\\
{\bf Remark}:
Note from \eqref{refequid1} and \eqref{solpart} that the barycenter of $\zeta_i(s)$ is fixed, in the sense that
\begin{equation}\label{barycenter}
\frac{\zeta_1(s)+ \dots +\zeta_k(s)}k= \frac{\bar\zeta_1(s)+ \dots +\bar\zeta_k(s)}k+\zeta_0=\zeta_0,\;\;\forall s\ge -\log T(0).
\end{equation}
Note that unlike in the one-dimensional case with a pure power nonlinearity treated in \cite{CZarxiv}, we are unable to prescribe the barycenter. Indeed, our equation \eqref{equ} is not invariant under the Lorentz transform.\\
{\bf Remark}: We are unable to say whether this solution has other characteristic points or not. In particular, we have been unable to find a solution with $\SS$ exactly equal to $\{0\}$. Nevertheless, let us remark that from the finite speed of propagation, we can prescribe more characteristic points, as follows:
\begin{coro}[Prescribing more characteristic points]\label{cormore} Let $J=\{1,...,n_0\}$ or $J=\N$ and for all $n\in J$, $r_n>0$, $T_n>0$ and $k_n \ge 2$
such that
\begin{equation}\label{necess}
r_n+T_n<r_{n+1}-T_{n+1}.
\end{equation}
Then, there exists a blow-up solution $u(r,t)$ of equation \eqref{gen}
with $\{x_n\;|\; n\in J\} \subset \SS$, $T(r_n)=T_n$ and for all $n\in I$,
\begin{equation*}
\left\|\vc{w_{x_n}(s)}{\ps w_{x_n}(s)} - \vc{\ds\sum_{i=1}^{k_n} (-1)^{i+1}\kappa(d_{i,n}(s))}0\right\|_{\H} \to 0\mbox{ as }s\to \infty,
\end{equation*}
with
\begin{equation*}
\forall i=1,\dots,k_n,\;\;d_{i,n}(s) = -\tanh \zeta_{i,n}(s),\;\;
\zeta_{i,n}(s) = \bar \zeta_i(s) + \zeta_{0,n}
\end{equation*}
for some $\zeta_{0,n}\in \R$, where $\bar \zeta_i(s)$ is defined in \eqref{solpart}.
\end{coro}
{\bf Remark}:  Again, we are unable to construct a solution with $\SS = \{r_n\;|\; n\in J\}$. When $N=1$, we may take $r_0\in \R$.

\bigskip

{\it Proof of Theorem \ref{mainth} and Corollary \ref{cormore}}: First, note that thanks to condition \eqref{necess} which asserts that the sections at $t=0$ of the backward light cones with vertices $(r_n, T_n)$ do not overlap, Corollary \ref{cormore} follows from Theorem \ref{mainth} by the finite speed of propagation. As for the proof of Theorem \ref{mainth}, we claim that it follows like in \cite{CZarxiv}, since the ingredients of that paper are available here, thanks to the adaptations we performed in the previous sections:\\
- {\bf the analysis of the ODE system \eqref{eqz}}: let us emphasize the fact that we still encounter this sytem in our setting. Indeed, that system appears as a projection on the null modes of the linearization of equation \eqref{eqw} around the sum of decoupled solitons, and, as we said in Part 2 page \pageref{pagepart2}, that equation differs from the pure power case, only with exponentially small terms (see \eqref{small}), which are absorbed in the $O(\frac 1{s^{1+\eta}})$ in \eqref{eqz};\\
- {\bf a reduction to a finite dimensional problem}: this is done thanks to the analysis of the dynamics of the equation in similarity variables \eqref{eqw} around the sum of decoupled solitons, which we did already for the proof  of the isolatedness of characteristic points (see (i) of Theorem \ref{thgeo}; see \cite{MZisol10} for the analysis in one space dimension);\\
- {\bf a topolgical argument to solve the  finite dimensional problem}: this argument is based on a different formulation of Brouwer's Theorem. It is independent of the equation.\\
%
%
%
%
%
%
Note however that one argument of \cite{CZarxiv} does not work here: the argument that allows us to prescribe the barycenter of $\zeta_i(s)$. Indeed, that argument uses the invariance of the pure power wave equation under the Lorentz transform, which is no longer the case for equation \eqref{equ}.\Box

\bigskip

Let us give in the following a criterion about the existence of characteristic points:
\begin{theor}[Existence and generic stability of characteristic points]
\label{thexis}
$ $\\
(i) {\bf (Existence)} Let $0<a_1<a_2$ be two non-characteristic points such that
\[
w_{a_i}(s) \to \theta(a_i)\kappa(d_{a_i},\cdot)\mbox{ as }s\to \infty\mbox{ with }\theta(a_1)\theta(a_2)=-1
\]
for some $d_{a_i}$ in $(-1,1)$, in the sense \eqref{profile}. Then, there exists a characteristic point $c\in (a_1,a_2)$.\\
(ii) {\bf (Stability)} There exists $\epsilon_0>0$ such that if $\|(\tilde U_0,\tilde U_1)- (U_0, U_1)\|_{\h1\times \l2(\R^N)}\le \epsilon_0$, then, $\tilde u(r,t)$ the solution of equation \eqref{equ} with initial data $(\tilde u_0,\tilde u_1)(r)=(\tilde U_0,\tilde U_1)(x)$ if $r=|x|$ blows up and has a characteristic point $\tilde c\in [a_1,a_2]$.
\end{theor}
{\bf Remark}: This statement (valid for $N\ge2$) is different from the original one (Theorem 2 in \cite{MZajm10}) by two natural small facts:
we take positive points $a_1$ and $a_2$ in (i), and
we use the multi-dimensional norm in (ii) (of course, from the finite speed of propagation, it is enough to take a localized norm instead). When  $N=1$, we don't need the restriction $a_1>0$.\\
{\bf Remark}: If one needs a quick argument for the existence of a blow-up solution for equation \eqref{equ} with a characteristic point, then this theorem allows us to avoid the heavy machinery of \cite{MZisol10}, namely the linearization of equation \eqref{eqw} around the sum of decoupled solitons. Indeed, we  have a more elementary argument, based on the knowledge of the blow-up behavior at a non-characteristic point on the one hand, and on (i) of this theorem on the one hand. However, such an argument uses  the fact that solutions of the ODE \eqref{ODE} associated to \eqref{equ} are also solution to \eqref{equ}
and this is possible only if $g$ defined in \eqref{gen} does not depend on $|x|$. For the statement with no perturbations, see Proposition 3 page 362 in \cite{MZbsm11}. For a further justification, see the Proof of Theorem \ref{thexis} below.

\medskip

{\it Proof of
Theorem \ref{thexis}}: As in \cite{MZbsm11}, there is no difficulty in adapting to the present context the proof of Theorem 2 of \cite{MZajm10} given in Section 2 of that paper, except may be for some natural extensions to the radial case.
Concerning the second remark following Theorem \ref{thexis},
the only delicate point is to find initial data $(u_0,u_1)$ satisfying the
hypothesis of (i) in Theorem \ref{thexis}. If $g$ does not depend on $|x|$, then, any solution of the ODE
\begin{equation}\label{ODE}
 U'' =|U|^{p-1}U+f(U)+g(t,0,U'),
\end{equation}
   is also a solution of the PDE \eqref{equ}, and it is enough to take  initial data $(u_0,u_1)$ with large plateaus of
 opposite signs. If $g$ does depend on $|x|$, then this simple idea breaks down, and the existence of initial data with characteristic points holds thanks to Theorem \ref{mainth}.
\Box

\bigskip

We also have the following result which relates the existence of characteristic points to the sign-change of the solution:
\begin{theor}
{\bf (Non-existence of characteristic points if the sign is constant)}
 Consider $u(r,t)$ a blow-up solution of \eqref{equ} such that $u(r,t)\ge 0$
  for all $r\in (a_0,b_0)$ and $t_0\le t< T(r)$ for some real $0\le a_0<b_0$ and $t_0\ge 0$. Then, $(a_0, b_0)\subset \RR$.
\end{theor}
{\bf Remark}: When $N=1$, we don't need the restriction $a_0\ge 0$.\\
{\it Proof}: This result follows from Theorem \ref{new} above
exactly as in one space dimension with no perturbations (i.e. for
equation (\ref{equ}) with $(f,g)\equiv(0,0)$). See the proof of
Theorem 4 given in Section 4.2 in \cite{MZajm10}.\Box

\appendix

\section{A local energy estimate for perturbations of the semilinear wave equation}\label{appenergy}
Let us consider the following perturbation of equation \eqref{par}:
\begin{equation}\label{genl}
\partial^2_{t} U =\Delta U+|U|^{p-1}U+
\lambda^{\frac{2p}{p-1}}f\Big(\lambda^{\frac{-2}{p-1}}U\Big) +
\lambda^{\frac{2p}{p-1}}g\Big(\lambda |x|,\lambda
t,\lambda^{-\frac{p+1}{p-1}}\nabla U\cdot \frac x{|x|},\lambda^{-\frac{p+1}{p-1}}\partial_tU\Big),
\end{equation}
where $\lambda>0$ and $f$ and $g$ satisfy ($H_f$) and ($H_g$). The equation \eqref{genl} is derived from equation \eqref{gen} through the dilation
\[
\lambda \mapsto U_\lambda(x,t) = \lambda^{\frac 2{p-1}}U(\lambda x, \lambda t).
\]
Using the technique of Shatah and Struwe \cite{SSnyu98} and introducing
\begin{eqnarray}
 \E(U(t)) &=&\int_{|x|<1-t}\left[\frac{(\partial_t U(x,t))^2}2+\frac{(\nabla U(x,t))^2}2 -\frac{|U(x,t)|^{p+1}}{p+1}-\lambda^{\frac{2p+2}{p-1}}F(\lambda^{\frac{-2}{p-1}}U(x,t))\right]{\mathrm{d}}x,\qquad\qquad
\end{eqnarray}
 where $F$ is defined in \eqref{F}.
We obtain the following local energy estimate
\begin{lem}[A local energy estimate for perturbations of equation \eqref{par}]
For all $t\in [0,1)$, we have
\begin{eqnarray*}
 \E(U(t))&\le& C\E(U(0))  +
C \int_0^t \int_{B_{s}} |U(\sigma ,s)|^{p+1}{\mathrm{d}}\sigma {\mathrm{d}}s
+C\lambda \int_0^t \int_{|x|<1-s} |U(x,s)|^{p+1}
{\mathrm{d}}x{\mathrm{d}}s+C\lambda^{\frac{2}{p-1}}.
\end{eqnarray*}
where the lateral boundary is
\begin{equation*}
B_{t_0}= \{(x,t)\;|\;0\le t \le t_0,\;|x|=1-t\}.
 \end{equation*}
\end{lem}
{\it Proof}: Classical  calculation  implies
\begin{eqnarray*}
\frac{d}{dt} \E(U(t))
&=&\underbrace{-\frac12 \int_{B_{t}}\Big(\partial_t U-\frac x{|x|}.\nabla U|\Big)^2 {\mathrm{d}}\sigma}_{L_1(s)} \underbrace{-\frac12\int_{B_{t}}\Big(|\nabla U|^2-(\frac x{|x|}.\nabla U)^2\Big){\mathrm{d}}\sigma}_{L_2(s)} \\
&&+\int_{B_{t}}\Big(\frac{|U|^{p+1}}{p+1}+
\lambda^{\frac{2p+2}{p-1}}F(\lambda^{\frac{-2}{p-1}}U(x,t))\Big){\mathrm{d}}\sigma\\
&&-\lambda^{\frac{2p}{p-1}}\int_{|x|<1-t}
g\Big(\lambda |x|,\lambda
t,\lambda^{-\frac{p+1}{p-1}}\nabla U\cdot \frac x{|x|},\lambda^{-\frac{p+1}{p-1}}\partial_tU\Big)\partial_tU {\mathrm{d}}x.
\end{eqnarray*}
Since $|\frac {x}{|x|}.\nabla U|\le |\nabla U|$, we can say that the term $L_2(s)$ is negative.
By combining  the estimate $|{F(x)}|\le  C( |x|+|x|^{p+1})$, the assumption  $(H_g)$
and the fact that the terms $L_1(s)$ and $L_2(s)$  are negative, we conclude that, for all $\lambda\in (0,1]$
\begin{eqnarray*}
\frac{d}{dt} \E(U(t))&\le&
C \int_{B_{t}}\Big(\lambda^{\frac{2p}{p-1}} |U|+ |U|^{p+1}\Big){\mathrm{d}}\sigma
+ C \lambda \int_{|x|<1-t}\Big( \lambda^{\frac{p+1}{p-1}}|\partial_tU|+
|\nabla U|^2+|\partial_tU|^2\Big){\mathrm{d}}x.
\end{eqnarray*}
\begin{eqnarray*}
\frac{d}{dt} \E(U(t))&\le& C\lambda  \E(U(t))+
C \int_{B_{t}}\Big(\lambda^{\frac{2p}{p-1}} |U|+ |U|^{p+1}\Big){\mathrm{d}}\sigma\\
&&+C\lambda \int_{|x|<1-t}
\Big( \lambda^{\frac{p+1}{p-1}}|\partial_tU|+\lambda^{\frac{2p}{p-1}} |U|+ |U|^{p+1}\Big)
{\mathrm{d}}x
\end{eqnarray*}
So
\begin{eqnarray*}
\frac{d}{dt} \E(U(t))\le C\lambda  \E(U(t))+
C \int_{B_{t}} |U|^{p+1}{\mathrm{d}}\sigma
+C\lambda \int_{|x|<1-t}
 |U|^{p+1}
{\mathrm{d}}x+C\lambda^{\frac{p+1}{p-1}}
\end{eqnarray*}
Then
\begin{eqnarray*}
 \E(U(t))&\le& \E(U(0))e^{C\lambda t}  +
C \int_0^te^{C\lambda (t-s) } \int_{B_{s}} |U|^{p+1}{\mathrm{d}}\sigma {\mathrm{d}}s\\
&&+C\lambda \int_0^te^{C\lambda (t-s)} \int_{|x|<1-s} |U|^{p+1}
{\mathrm{d}}x{\mathrm{d}}s+C\lambda^{\frac{p+1}{p-1}} \int_0^te^{C\lambda (t-s)}  {\mathrm{d}}s
\end{eqnarray*}
so
\begin{eqnarray*}
 \E(U(t))&\le& \E(U(0))e^{C\lambda t}  +
C \int_0^te^{C\lambda (t-s) } \int_{B_{s}} |U|^{p+1}{\mathrm{d}}\sigma {\mathrm{d}}s\\
&&+C\lambda \int_0^te^{C\lambda (t-s)} \int_{|x|<1-s} |U|^{p+1}
{\mathrm{d}}x{\mathrm{d}}s+C\lambda^{\frac{2}{p-1}} e^{C\lambda t}
\end{eqnarray*}

............. \Box
\def\cprime{$'$}


\noindent{\bf Address}:\\
Universit\'e de Tunis El-Manar, Facult\'e des Sciences de Tunis, D\'epartement de math\'ematiques, Campus Universitaire 1060,
 Tunis, Tunisia.\\
\vspace{-7mm}
\begin{verbatim}
e-mail: ma.hamza@fst.rnu.tn
\end{verbatim}
Universit\'e Paris 13, Sorbonne Paris Cit\'e, LAGA, CNRS (UMR 7539),
99 avenue J.B. Cl\'ement, 93430 Villetaneuse, France.\\
\vspace{-7mm}
\begin{verbatim}
e-mail: Hatem.Zaag@univ-paris13.fr
\end{verbatim}

\end{document}